\newcommand{\cen}{\begin{displaymath}}
\newcommand{\ter}{\end{displaymath}}
\newcommand{\hz}{\vspace{0.5cm}}
\newcommand{\kla}{\left ( }
\newcommand{\mer}{\right ) }
\newcommand{\ruei}{RU\!M\!D_n^1(}
\newcommand{\rup}{RU\!M\!D_n^p(}
\newcommand{\ruq}{RU\!M\!D_n^q(}
\newcommand{\rur}{RU\!M\!D_n^r(}
\newcommand{\ruz}{RU\!M\!D_n^2(}
\newcommand{\rue}{RU\!M\!D_n (}
\newcommand{\ruee}{RU\!M\!D_{2n} (}
\newcommand{\vare}{\varepsilon}
\newtheorem{lemma}{Lemma}[section]

\newtheorem{theorem}[lemma]{Theorem}
\newtheorem{cor}[lemma]{Corollary}
\newcommand{\ew}{{\rm  I\! E}}
\newcommand{\nz}{{\rm  I\! N}}

\newcommand{\kz}{{\rm  I\! K}}
\newcommand{\pl}{\hspace{.1cm}}
\newcommand{\pll}{\hspace{.3cm}}
\newcommand{\pla}{\hspace{1.5cm}}
\newcommand{\omet}{\Omega_n }
\newcommand{\omett}{\Omega_{2n} }
\newcommand{\omep}{I\! D_n }
\documentstyle{article}
\begin{document}

\title{Lower estimates of random unconditional constants of Walsh-Paley
martingales with values in Banach spaces}

\author{Stefan Geiss \thanks{The author is supported by the DFG
(Ko 962/3-1).} \\
Mathematisches Institut der Friedrich Schiller Universit\"at Jena\\
UHH 17.OG, D O-6900 Jena, Germany} 

\maketitle

{\bf Abstract.} For a Banach space $X$ we define $\rue X)$ to be the infimum 
of all $c>0$ such that 
\cen \kla AV_{\vare_k =\pm 1} \| \sum_1^n \vare_k (M_k - M_{k-1} )\|_{L_2^X}^2
\mer^{1/2} \le c \| M_n \|_{L_2^X} \ter
holds for all Walsh-Paley martingales $\{ M_k \}_0^n \subset L_2^X$ with
$M_0 =0$. We relate the asymptotic behaviour of the sequence 
$\{ \rue X)\}_{n=1}^{\infty}$ to geomertrical properties of the Banach space
$X$ such as K-convexity and superreflexivity.\hz

1980 Mathematics Subject Classification (1985 Revision): 46B10, 46B20

\setcounter{section}{-1}

\section{Introduction}
\setcounter{lemma}{0}

A Banach space $X$ is said to be an UMD-space if for all $1<p<\infty$ there is 
a constant $c_p = c_p (X)>0$ such that

\cen \sup_{\vare_k =\pm 1}\| \sum_1^n \vare_k (M_k -M_{k-1} ) \|_{L_p^X}
\le c_p \| \sum_1^n (M_k - M_{k-1} )\|_{L_p^X} \ter

for all $n=1,2,...$ and all martingales $\{ M_k \}_0^n \subset L_p^X$ with 
values in $X$. It turns out that
this definition is equivalent to the modified one if we replace "for all
$1<p<\infty$" by "for some $1<p<\infty$", and "for all martingales" by "for all
Walsh-Paley-martingales" (see \cite{Bu} for a survey). Motivated by these
definitions we investigate Banach spaces $X$ by means of the sequences
$\{ \rue X)\}_{n=1}^{\infty}$ whereas $\rue X) := \inf c$ such that

\cen \kla AV_{\vare_k =\pm 1} \| \sum_1^n \vare_k (M_k - M_{k-1} )\|_{L_2^X}^2
\mer^{1/2} \le c \| M_n \|_{L_2^X} \ter

holds for all Walsh-Paley martingales $\{ M_k \}_0^n \subset L_2^X$ with the
starting point $M_0 =0$. "$RU\!M\!D$" stands for "random unconditional
constants of martingale differences". We consider "random" unconditional
constants instead of the usual one, where $\sup_{\vare =\pm 1}$ is taken in 
place of $AV_{\vare =\pm 1}$, since they naturally appear in the lower 
estimates we are interested in. These lower estimates are of course lower 
estimates for the non-random case, too. The paper is organized in the 
following way. Using a technique of Maurey we show that the exponent 2 in the 
definition of $\rue X)$ can be replaced by any $1<p<\infty$ (see Theorem 2.4). 
Then we observe (see Theorem 3.5)

\cen X \pll \mbox{is not K-convex} \hspace{1cm} \Longleftrightarrow 
\hspace{1cm} \rue x) \asymp n. \ter
 
In the case of superreflexive Banach spaces this turns into

\cen X \pll \mbox{is not superreflexive} \hspace{1cm} \Longrightarrow 
\hspace{1cm} \rue X) \succeq n^{1/2}, \ter

and, under the assumption X is of type 2,

\cen X \pll \mbox{is not superreflexive} \hspace{1cm} \Longleftrightarrow 
\hspace{1cm} \rue x) \asymp n^{1/2} \ter

(see Theorems 4.3 and 4.4). Using an example due to Bourgain we see that the
type 2 condition is necessary. In fact, for all $1<p<2<q<\infty$ there is a
superreflexive Banach space $X$ of type $p$ and cotype $q$ such that
$\rue X) \asymp n^{\frac{1}{p} - \frac{1}{q}}$ (see Corollary 5.4).
According to a result of James a non-superreflexive Banach spaces $X$ is 
characterized by the existence of large "James-trees" in the unit ball $B_X$ 
of $X$. We can identify these trees with Walsh-Paley martingales 
$\{ M_k \}_0^n$ which only take values in the unit ball $B_X$ and which satisfy
$\inf_{\omega} \| M_k (\omega ) - M_{k-1} (\omega )\| \ge \theta$ for
some fixed 
$0<\theta <1$. In this way we can additionally show
that the martingales, which give the lower estimates of our random 
unconditional constants, are even James trees (see Theorems 3.5(2) and 4.3).

 \section{Preliminaries}
\setcounter{lemma}{0} 

The standard notation of the Banach space theory is used
(cf.\cite{L-T}). Throughout this paper $\kz$ stands for the real or complex 
scalars. $B_X$ is the closed unit ball of the Banach space $X$, ${\cal L}
(X,Y)$ 
is the space of all linear and continuous operators from a Banach space $X$ 
into a Banach space $Y$ equipped with the usual operator norm.
We consider martingales over the probability space 
$[\omet,\mu_n ]$ which is given by $\omet := \{ \omega =(\omega_1
,...,\omega_n ) \in \{ -1,1\} ^n \}$ and 
$\mu_n (\omega ) :=\frac{1}{2^n}$ for all $\omega \in \omet$.  The minimal 
$\sigma$-algebras ${\cal F} _k $, such that the coordinate functioinals
$\omega =(\omega_1 ,...,\omega_n ) \rightarrow \omega_i \in \kz$
are measurable for $i=1,...,k$, and ${\cal F}_0 :=\{\emptyset ,\omet \}$ 
form a natural filtration $\{ {\cal F} _k\}$ on $\omet$. 
A martingale $\{ M_k \}_0^n$ with values in a Banach space
$X$ over $[\omet ,\mu_n ]$ with respect to this filtration $\{ {\cal F} \}_0^n$
is called Walsh-Paley martingale.  As usual we put $dM_0 := M_0,
dM_k := M_k -M_{k-1}$ for $k\ge 1$ and $M_k^*(\omega ) = \sup_
{0\le l\le k} \| M_l(\omega )\|$.  Given a function $M \in L_p^X
(\omet )$ we can set $M_k :=\ew (M|{\cal F }_k )$ for
$k=0,...,n$. Consequently, for each $M \in L_p^X (\omet )$ there
is a unique Walsh-Paley martingale $\{ M_k \}_0^n$ with $M_n
=M$. In this paper we consider a further probability space
$[\omep ,P_n ]$ with $\omep =
\{ \vare =(\vare_1 ,...,\vare_n ) \in \{ -1,1\} ^n \}$ and
$P_n (\vare )= \frac{1}{2^n}$ for all $\vare \in \omep$.
$\ew_{\vare, \omega }$ means that we take the expectation with
respect to the product measure $P_n \times \mu_n$.  To estimate the
random unconditional constants of Walsh-Paley martingales from
above we use the notion of the type.  For $1\le p\le 2$ an
operator $T \in {\cal L }(X,Y)$ is of type $p$ if

\cen \left( \ew _{\vare } \| \sum _k T \vare_k x_k \| ^2 \right)
^{1/2} \le c \left( \sum_k \| x_k \|^p \right)^{1/p} \ter

for some constant $c>0$ and all finite sequences $\{x_k \} \subset X$. The
infimum of all possible constants $c>0$ is denoted by $T_p (T)$.
Considering the above inequality for sequences $\{ x_k \}_{k=1}^n \subset X$
of a fixed length $n$ only we obtain the corresponding constant $T_p^n (T)$
which can be defined for each operator $T\in {\cal L }(X,Y)$.
In the case $T=I_X$ is the identity of a Banach space $X$ we write
$T_p (X)$ and $T_p^n (X)$ instead of $T_p (I_X)$ and $T_p^n (I_X)$, and
say "$X$ is of type p" in place of "$I_X$ is of type $p$" (see \cite{T-J} for 
more information).
\pagebreak

\section{Basic definition}
\setcounter{lemma}{0}

Let $T\in {\cal L} (X,Y)$ and $1 \le q<\infty $. Then $\ruq T) = \inf c$,
where the infimum is taken over all $c>0$ such that
\cen \left( \ew_{\varepsilon ,\omega } \| \sum_1^n \varepsilon _k TdM_k
(\omega )\| ^q \right)^{1/q} \le c \left( \ew _{\omega } \| M_n (\omega )\|^q
\right) ^{1/q}\ter

holds for all Walsh-Paley martingales $\{ M_k\}_0^n$ with values in $X$ and
$M_0 =0$. Especially, we set \linebreak
$\ruq X) := \ruq I_X )$ for a Banach space $X$ with the identity $I_X$.
It is clear that $\ruq T) \le 2n \| T \|$ since
\cen \| \sum_1^n \vare_k TdM_k \|_{L_q^X} \le \sum_1^n \| T\| 
\| dM_k \|_{L_q^X} \le  \| T\| \sum_1^n \kla \| M_k\|_{L_q^X} + 
\| M_{k-1}\|_{L_q^X} \mer \le  2 n \| T\| \| M_n \|_{L_q^X}. \ter
In the case $X$ is an UMD-space we have $\sup_n \ruq X) < \infty$ whenever 
$1<q< \infty$ (the converse seems to be open). $q=1$ yields a"singularity" 
since $\ruei X) \asymp \ruei \kz ) \asymp n$ for any Banach space $X$ (see 
Corollary 5.2) therefore we restrict our consideration on $1<q< \infty$. Here 
we show that the quantities $\ruq T)$ are equivalent for $1<q<\infty$. In 
\cite{Ga}(Thm.4.1) it is stated that $\sup_n \ruq X)<\infty$ iff $\sup_n 
\rur X)<\infty$ for all $1<q,r<\infty$. Using Lemma 2.2, which slightly extends
\cite{M}(Thm.II.1), we prove a more precise result in Theorem 2.4.
\hz

Let us start with a general martingale transform. Assuming
$T_1 ,...,T_n \in {\cal L}(X,Y)$ we define

\cen \phi = \phi (T_1 ,...,T_n ) : L_0^X (\omet ) \longrightarrow
L_0^Y (\omet ) \pll \mbox{by} \pll
\phi (M)(\omega ):=\sum_1^n T_k dM_k (\omega ),\ter

where $M_k = \ew (M|{\cal F}_k )$. The following duality is standard.\hz

\begin{lemma} Let $1<p<\infty$,\pl $T_1,...,T_n \in {\cal L}(X,Y)$ and
$\phi =\phi (T_1 ,...,T_n ) : L_p^X \longrightarrow L_p^Y $. Then

\cen \phi '(F) = \sum_1^n T_k 'dF_k \pll \mbox{for all} \pll
F \in L_{p'}^{Y'}.\ter

\end{lemma} \hz

$Proof.$ Using the known formula

\cen <M,M'> = \sum_0^n <dM_k ,{dM'}_k >
\pll \kla M \in L_s^Z (\omet ),M'\in L_{s'}^{Z'} (\omet ),1<s<\infty \mer,
\ter

for $M\in L_p^X (\omet )$ and $F\in L_{p'}^{Y'} (\omet )$ we obtain

\begin{eqnarray*}
<M,\phi 'F> & = & <\sum_1^n T_k dM_k ,F> = \sum_1^n <T_k dM_k ,dF_k > \\
& = & \sum_1^n <dM_k ,T_k 'dF_k > = <M,\sum_1^n T_k 'dF_k >.\pll \Box
\end{eqnarray*} \hz

Now we recall \cite{M}(Thm.II.1) in a more general form. Although the proof 
is the same we repeat some of the details for the convenience of the reader.\hz

\begin{lemma} Let $1<p<r<\infty$, $T_1 ,...,T_n \in {\cal L}(X,Y)$ and
$\phi = \phi (T_1 ,...,T_n )$. Then
\cen \| \phi :L_p^X \rightarrow L_p^Y \| \le \frac{6r^2}{(p-1)(r-1)}
\| \phi :L_r^X \rightarrow L_r^Y \| .\ter \end{lemma} \hz

$Proof.$ We define $1<q< \infty$ and $0<\alpha <1$ with $\frac{1}{p} =
\frac{1}{r} +\frac{1}{q} $ and $\frac{p}{r} = 1- \alpha $ and obtain
$\alpha q=(1-\alpha )r=p$. Let $\{ M_k \} _0^n$ be a Walsh-Paley martingale
in $X$ with $dM_1 \neq 0$. We set

\cen ^*\! M_k (\omega ) := M_{k-1}^* (\omega ) +
sup_{0\le l \le k} \| dM_l (\omega )\| \hspace{1cm}for\hspace{1cm}
k\ge 1 \ter

and obtain an ${\cal F} _{k-1}^{\omega }$ -measurable random variable with

\cen 0< {^*\! M}_1(\omega )\le ... \le {^*\! M}_n(\omega ) \pll \mbox{and} \pll
M_k^*(\omega ) \le {^*\! M}_k(\omega ) \le 3 M_k^*(\omega ). \ter

Using \cite{M}(L.II.B) and Doob's inequality we obtain

\begin{eqnarray*}
\| \phi (M) \|_p & = &
\left( \ew _{\omega }\| \sum_1^n
\frac{T_k dM_k (\omega )}{^*\! M_k^{\alpha } (\omega)} {^*\! M}_k^{\alpha }
(\omega ) \| ^p \right) ^{1/p} \\
& \le & 2 \left( \ew _{\omega } \sup_{1\le k\le n} \| \sum _1^k
\frac{T_l dM_l (\omega )}{^*\! M_l^{\alpha }(\omega )} \| ^p \hspace{0.2cm}
{^*\! M_n^{\alpha p}}(\omega) \right) ^{1/p} \\
& \le & 2 \left( \ew _{\omega } \sup_{1\le k\le n} \| \sum _1^k
\frac{T_l dM_l (\omega )}{^*\! M_l^{\alpha }(\omega )} \| ^r \right) ^{1/r}
\left( {\ew _{\omega }} ^*\! M_n^{\alpha q} (\omega ) \right) ^{1/q} \\
& \le & \frac{2r}{r-1} \left( \ew _{\omega } \| \sum _1^n
\frac{T_k dM_k (\omega )}{^*\! M_k^{\alpha }(\omega )} \| ^r \right) ^{1/r}
\left( {\ew _{\omega }} ^*\! M_n^p (\omega ) \right) ^{1/q} \\
& \le & \frac{2r}{r-1} \| \phi \|_r \left( \ew _{\omega } \| \sum _1^n
\frac{dM_k (\omega )}{^*\! M_k^{\alpha }(\omega )} \| ^r \right) ^{1/r}
\left( {\ew _{\omega }} ^*\! M_n^p (\omega ) \right) ^{1/q}. \\
\end{eqnarray*}

Applying \cite{M}(L.II.A) in the situation  $ \| \sum _1^l dM_i (\omega )
\| \le \left( ^*\! M_l (\omega )^{\alpha } \right) ^{1/ {\alpha }} $ yields

\cen \| \sum _1^k \frac{dM_l (\omega )}{^*\! M_l^{\alpha }(\omega )} \| \le
\frac{1/{\alpha }}{1/{\alpha } -1} {^*\! M_k(\omega)}^
{{\alpha }(1/{\alpha } -1)} \le \frac{r}{p} {^*\! M_n(\omega )}^{p/r} \ter

such that

\begin{eqnarray*}
\| \phi (M)\|_p & \le & \frac{2r^2}{(r-1)p} \| \phi \|_r \kla \ew_{\omega }
^* M_n (\omega )^p \mer^{1/p}
\le \frac{6r^2}{(r-1)p} \| \phi \|_r \kla \ew_{\omega } M_n^* (\omega )^p
\mer^{1/p} \\
& \le & \frac{6r^2}{(r-1)(p-1)} \| \phi \|_r \| M_n \|_p. \pll \Box
\end{eqnarray*} \hz

We deduce \hz

\begin{lemma} Let $T_1 ,...,T_n :X\rightarrow L_1^Y (\omep )$ and
$\phi := \phi (T_1 ,...,T_n )$. For $i=1,...,n$ assume that \linebreak 
$T_i (X) \subseteq \{ f:f=\sum_1^n \vare_k y_k, \hspace{0.2cm} y_k \in Y \}$. 
Then
\cen \frac{1}{c_{pr}} \| \phi : L_p^X \rightarrow L_p^{L_p^Y} \| \le
\| \phi : L_r^X \rightarrow L_r^{L_r^Y} \| \le c_{pr}
\| \phi : L_p^X \rightarrow L_p^{L_p^Y} \| \ter
for $1<p<r<\infty$, where the constant $c_{pr} >0$ is independent from
$X$,$Y$,$(T_1 ,...,T_n )$ and $n$. \end{lemma} \hz

$Proof.$ The left-hand inequality follows from Lemma 2.2 and
\cen  \| \phi : L_p^X \rightarrow L_p^{L_p^Y} \| \le \frac{6r^2}{(p-1)(r-1)}
\| \phi : L_r^X \rightarrow L_r^{L_p^Y} \| \le  \frac{6r^2}{(p-1)(r-1)}
\| \phi : L_r^X \rightarrow L_r^{L_r^Y} \| .\ter
The right-hand inequality is a consequence of Lemma 2.1, Lemma 2.2 and
\begin{eqnarray*}
\| \phi : L_r^X \rightarrow L_r^{L_r^Y} \| & = &
\| \phi ': L_{r'}^{L_{r'}^{Y'}} \rightarrow L_{r'}^{X'} \| \\
& \le & \frac{{6r'}^2}{(p'-1)(r'-1)}
\| \phi ': L_{p'}^{L_{r'}^{Y'}} \rightarrow L_{p'}^{X'} \| \\
& = & \frac{{6r'}^2}{(p'-1)(r'-1)}
\| \phi : L_p^X \rightarrow L_p^{L_r^Y} \| \\
& \le & \frac{{6r'}^2}{(p'-1)(r'-1)} K_{rp}
\| \phi : L_p^X \rightarrow L_p^{L_p^Y} \|,
\end{eqnarray*} 
where we use Kahane's inequality (cf. \cite{L-T} (II.1.e.13)) in the
last step.\pll $\Box$ \hz

If we apply Lemma 2.3 in the situation $T_k x:= \vare _k Tx$ and exploit
\cen \rup T) \le \| \phi (T_1 ,...,T_n ):L_p^X \rightarrow L_p^{L_p^Y} \|
\le 2 \rup T) \ter
then we arrive at \hz

\begin{theorem} Let $1<p<r<\infty$ and $T \in {\cal L }(X,Y)$. Then
\cen \frac{1}{c_{pr}} \rup T) \le \rur T) \le c_{pr} \rup T), \ter
where the constant $c_{pr} >0$ is independent from $X$,$Y$,$T$ and $n$.
\end{theorem} \hz

The above consideration justifies
\cen \rue T):=\ruz T) \pll \mbox{for} \pll T \in {\cal L }(X,Y) \ter
and $\rue X):= \rue I_X )$ for a Banach space $X$.\\
\hz

\section{K--convexity}
\setcounter{lemma}{0}

We show that $\rue X) \asymp n$ if and only if $X$ is not K-convex, that is,
if and only if $X$ uniformly contains $l_1^n$. To do this some additional 
notation is required. For $x_1,...,x_n\in X$ we set
\cen |x_1\wedge ...\wedge x_n|_X:=\sup \{ |det(\langle x_i,a_j\rangle )
_{i,j=1}^n|:a_1 ,...,a_n \in B_{X'} \}.\ter
Furthermore, for fixed $n$ we define the bijection
\cen i: \{ -1,1\} ^n \rightarrow \{1,...,2^n \} \ter as
\cen i(\omega )=i(\omega_1,...\omega_n ):=1+ \frac{1-\omega_n}{2} +
\frac{1-\omega_{n-1}}{2} 2 +...+ \frac{1-\omega_1}{2} 2^{n-1}  \ter
and the corresponding sets $I_0 := \{ 1,...,2^n\},
I(\omega_1,...,\omega_n ):= \{ i(\omega_1 ,...,\omega_n )\}$,
\cen I(\omega_1,...,\omega_k ):=
\{ i(\omega_1,...,\omega_n): \omega_{k+1} =\pm 1,...,\omega _n =\pm 1\}
\hspace{1cm} \mbox{for} \hspace{1cm} k=1,...,n-1.\ter
It is clear that
\cen I(\omega_1,...,\omega_{k-1})=I(\omega_1,..., \omega_{k-1},1)
\cup I(\omega_1,...,\omega_{k-1},-1) \hspace{1cm} \mbox{and} \hspace{1cm}
I_0 = I(1) \cup I(-1). \ter

Our first lemma is technical. \hz

\begin{lemma} Let $\{ M_k\} _0^n$ be a Walsh-Paley-martingale in $X$ and let
$x_k := M_n(i^{-1} (k))\in X$ for \linebreak
$k=1,...,2^n$. Then, for all $\omega \in
\{ -1,1\} ^n$ and $1\le k\le n$ , there exist natural numbers \linebreak
$1\le r_0\le s_0<r_1 \le s_1 <...<r_k \le s_k \le 2^n$
with
\cen |M_0(\omega )\wedge ...\wedge M_k(\omega )|=\frac{1}{2^k} \left | \frac
{x_{r_0}+...+x_{s_0}}{s_0 -r_0 +1} \wedge ...\wedge \frac
{x_{r_k}+...+x_{s_k}}{s_k -r_k +1}\right |  \ter \end{lemma} \hz

$Proof.$  Let us fix $\omega \in \{ -1,1\} ^n.$ Since
$M_l(\omega _1,...,\omega _l) = \frac{1}{2^{n-l}} \sum_
{i\in I(\omega _1,...,\omega _l)} x_i$
we have for $l=0,...,n-1,$

\cen M_l(\omega _1,...,\omega _l) - \frac{1}{2} M_{l+1}
(\omega _1,...,\omega _{l+1}) =\frac{1}{2^{n-l}} \sum_
{I(\omega _1,...,\omega _l ,-\omega _{l+1})} x_i = 
\frac{1}{2\# I(\omega _1,...,\omega _l,
-\omega _{l+1})}\sum_{I(\omega _1,...,\omega _l ,
-\omega _{l+1})} x_i \ter

for $l=0,...,n-1$. It is clear that $I(-\omega _1),
I(\omega _1,-\omega _2),
I(\omega _1,\omega _2,-\omega _3),...,
I(\omega _1,...,\omega _{k-1},-\omega _k),
I(\omega _1,...,\omega _k)$  are disjoint, such that we have

\begin{eqnarray*}
|M_0(\omega)\wedge ...\wedge M_k(\omega )| & = &
\left| \left( M_0(\omega)-\frac{1}{2} M_1(\omega) \right) \wedge ...
\wedge \left( M_{k-1}(\omega ) -\frac{1}{2} M_k(\omega ) \right)
\wedge M_k(\omega ) \right| \\ & = &
\frac{1}{2^k} \left | \frac{x_{r_0}+...+x_{s_0}}{s_0 -r_0 +1}
\wedge ...\wedge \frac{x_{r_k}+...+x_{s_k}}{s_k -r_k +1} \right |
\end{eqnarray*}

after some rearrangement. $\Box$ \hz

The second lemma, which is required, is a special case of \cite{G} (Thm.1.1).

\begin{lemma} Let $u \in {\cal L} (l_2^n ,X)$ and let $\{ e_1 ,...,e_n \}$
be the unit vector basis of $l_2^n$. Then
\cen |ue_1\wedge ...\wedge ue_n|_X \le \kla \frac{1}{n!} \mer^{1/2}
\pi_2(u)^n ,\ter
where $\pi_2(u)$ is the absolutely 2-summing norm of $u$.
\end{lemma} \hz

Now we apply Lemma 3.1 to a special Walsh-Paley-martingale $\{ M_k^1\}_0^n$
with values in $l_1^{{2^n}}$ whose differences $dM_k^1(\omega )$ are closely
related to a discret version of the Haar functions from $L_1 [0,1]$.  
For fixed $n$ this martingale is given by
\cen M_n^1(\omega _1,...,\omega _n) := e_{i(\omega _1,...,\omega _n)} 
\hspace{1cm} \mbox{and} \hspace{1cm}
M_k^1 := \ew (M_n^1 |{\cal F}_k ), \ter
where $\{ e_1,...,e_{2^n}\} $ stands for the unit vector basis of
$l_1^{{2^n}}.$ \hz

\begin{lemma} Let $n\ge 1$ be fixed. Then \\
(1) $\| M_k^1 (\omega )\| = \| dM_k^1 (\omega )\| = 1 $ \pll for
    $k=0,...n$ and all $\omega \in \omet$,\\
(2) $\inf_{\omega } \ew_{\vare } \| \sum_1^n \vare_k dM_k^1(\omega ) \|
    \ge \alpha n$ \pll for some $\alpha >0$ independent from $n$.
\end{lemma} \hz

$Proof.$ (1) is trivial. We consider (2). Lemma 3.1 implies

\cen \frac{1}{2^n} |f_1\wedge ...\wedge f_{n+1}|_{l_1^{n+1}} =
|M_0^1(\omega)\wedge ...\wedge M_n^1(\omega )|_{l_1^{2^n}} \ter

for all $\omega \in \{ -1,1\} ^n$, where $\{ f_1,...,f_{n+1}\} $ denotes
the unit vector basis of $l_1^{n+1}$. We can continue to

\begin{eqnarray*}
\frac{1}{2} |f_1\wedge ...\wedge f_{n+1}|_{l_1^{n+1}}^{1/n} & = &
|M_0^1(\omega)\wedge dM_1^1(\omega)\wedge ...\wedge dM_n^1(\omega )|
_{l_1^{2^n}}^{1/n} \\ & \le &
\left( (n+1)\| M_0^1(\omega) \| \hspace{0.3em} |dM_1^1(\omega)\wedge...
\wedge dM_n^1(\omega )| \right) ^{1/n} \\ & \le &
c |dM_1^1(\omega)\wedge...\wedge dM_n^1(\omega )| ^{1/n}.
\end{eqnarray*}

If we define the operator $u_{\omega } : l_2^n \longrightarrow l_1^{2^n}$ by
$u_{\omega } ((\xi_1 ,...,\xi_n )) := \sum_1^n \xi_i dM_i^1 (\omega )$ and use
Lemma 3.2, then we get

\cen |dM_1^1(\omega)\wedge...\wedge dM_n^1(\omega )| ^{1/n} \le
\left( \frac{1}{n!} \right) ^{1/2n} \pi _2 (u_{\omega }). \ter

Since the $l_1^{2^n}$ are uniformly of cotype 2 there is a constant $c_1 >0$, 
independent from $n$, such that

\cen \pi _2 (u_{\omega }) \le c_1 \ew_{\vare } \| \sum _1^n \vare _k
dM_k^1(\omega ) \| \ter

(see \cite{T-J},\cite{M-P}). Summerizing the above estimates yields

\cen |f_1\wedge ...\wedge f_{n+1}|_{l_1^{n+1}}^{1/n} \le 2ce^{1/2} n^{-1/2}
c_1 \ew_{\vare } \| \sum _1^n \varepsilon _k dM_k^1(\omega ) \| \le c_2
n^{-1/2} \ew_{\vare } \| \sum _1^n \varepsilon _k dM_k^1(\omega ) \|. \ter

The known estimate
$|f_1\wedge ...\wedge f_{n+1}|^{1/{n+1}} \ge \frac{1}{c_3}(n+1)^{1/2} $
concludes the proof (see, for instance, \cite{G}(Ex.2.7)).\\
\hspace*{\fill} $\Box$ \\

Finally, we need the trivial \hz

\begin{lemma} Let $T \in {\cal L }(X,Y)$. Then $\rue T) \le 2 n^{1/2}
T_2^n (T)$. \end{lemma} \hz

$Proof.$ Using the type 2 inequality for each $\omega \in \omet$ and
integrating yield for a martingale $\{ M_k \} _0^n$
 
\cen 
\left( \ew _{\vare ,\omega }\| \sum_1^n \vare_k T dM_k (\omega )\| ^2 
\right) ^{1/2} \le  
T_2^n(T) \left( \sum _1^n \| dM_k \|_{L_2^X}^2 \right) ^{1/2} \le 
2 n^{1/2} T_2^n (T) \| M_n \| _{L_2^X} .\Box
\ter

Now we can prove \hz

\begin{theorem} There exists an absolute constant $\alpha >0$ such that for
any Banach space $X$ the following assertions are equivalent.\\
(1) $X$ is not K-convex.\\
(2) For all $\theta >0$ and all $n=1,2,...$ there is a Walsh-Paley martingale
$\{ M_k\}_0^n$ with values in $B_X$,
\cen
\inf_{1\le k\le n} \inf_{\omega } \| dM_k (\omega )\| \ge 1-\theta
\pll \mbox{and} \pll
\inf_{\omega } \ew_{\vare } \| \sum_1^n \vare_k dM_k(\omega ) \| \ge \alpha n.
\ter
(3) $\rue X)\ge cn$ for $n=1,2,...$ and some constant $c=c(X)>0$.
\end{theorem} \hz

$Proof.$ Taking $\alpha >0$ from Lemma 3.3 the implication 
(1) $\Rightarrow$ (2) follows.
(2) $\Rightarrow$ (3) is trivial.\linebreak
(3) $\Rightarrow$ (1): Assumig $X$ to be K-convex the space $X$ must be of
type p for some $p>1$. Consequently, Lemma 3.4 implies

\cen \rue X)\le 2 n^{1/2} T_2^n(X) \le 2 n^{1/2} n^{1/p - 1/2}
T_p(X) \le 2 n^{1/p} T_p(X). \Box \ter

$Remark.$ One can also deduce (3) $\Rightarrow$ (1) from \cite{El} and
\cite{Pa} in a more direct way (we would obtain that $L_2^X (\{ -1,1\}^{\nz} )$
is not K-convex).

\section{Superreflexivity}
\setcounter{lemma}{0}

A Banach space $X$ is superreflexive if each Banach space, which is 
finitely representable in $X$, is reflexive. We will see that 
$\rue X)\ge cn^{1/2} $ whenever $X$ is not superreflexive and
that the exponent $\frac{1}{2}$ is the best possible in general. This improves
an observation of Aldous and Garling (proofs of \cite{Ga}(Thm.3.2) and 
\cite{Al}(Prop.2)) which says that $\rue X) \ge c n^{1/s}$ in the case $X$ is 
of cotype s ($2 \le s < \infty $) and not superreflexive. \hz

We make use of the summation operators

\cen \sigma _n : l_1^{2^n} \longrightarrow l_{\infty }^{2^n} \hspace{1cm}
and \hspace{1cm} \sigma : l_1 \longrightarrow l_{\infty }
\hspace{1cm} \mbox{with} \hspace{1cm} 
\{ \xi _k \}_k \longrightarrow \left\{ \sum _{l=1}^k \xi _l \right\} _k ,\ter
as well as of
\cen \Phi :C[0,1]' \longrightarrow l_{\infty } ([0,1]) \hspace{1cm} \mbox{with}
\hspace{1cm} \mu \longrightarrow \left\{ t \rightarrow \mu ([0,t]) \right\} .
\ter

The operators $\sigma_n$ are an important tool in our situation. Assuming $X$ 
to be not superreflexive, according to \cite{J1} for all $n=1,2,...$ there are 
factorizations $\sigma_n = B_n A_n $ with $A_n :l_1^{2^n} \rightarrow X$,
$B_n :X \rightarrow l_{\infty }^{2^n}$ and
$\sup _n \| A_n \| \| B_n \| \le 1+\theta\pl (\theta >0)$. It turns out that
the image-martingale $\{ M_k \}_0^n \subset L_2^X$ of $\{ M_k^1 \}$ 
($n$ is fixed, $\{ M_k^1 \}$ is defined in the previous section), which is 
given by $M_k (\omega ):=A_n M_k^1 (\omega ) \hspace{.3cm} (k=1,...,n)$, 
possesses a large random unconditional constant. To see this we set

\cen M_k^{\infty }(\omega ) := \sigma _n M_k^1(\omega ) \hspace{2cm}
\left( \omega \in \omet , k=0,...,n \right) \ter

and obtain a martingale $\{ M_k^{\infty }\} _0^n $ with values in
$l_{\infty }^{2^n}$. For $k=1,...,n$ it is easy to check that

\cen dM_k^{\infty }(\omega _1 ,..., \omega _k ) = \omega _k 2^{k-n-1}
(0,...,0,1,2,3,...,2^{n-k},2^{n-k} -1,...,3,2,1,0,0,...,0) \ter

where the block $(1,2,3,...,2^{n-k})$ is concentrated on
$I(\omega _1 ,..., \omega _{k-1} ,1)$ and the block $(2^{n-k} -1,...,3,2,1,0)$
is concentrated on $I(\omega _1 ,..., \omega _{k-1} ,-1)$, that is, the vectors
$|dM_k^{\infty} (\omega )|$ correspond to a discrete Schauder
system in $l_{\infty}^{2^n}$. Furthermore, we have \hz

\begin{lemma} Let $n \ge 2$ be a natural number and let $\{ e_i \}$ be the
standard basis of $l_1^{2^n}$. Then there exists a map
$e:\{ -1,1\} ^n \longrightarrow \{ e_1 ,...,e_{2^n} \} \subset l_1^{2^n}$
such that

\cen \mu_n \left\{ \omega : |<dM_k^{\infty }(\omega ),e(\omega )>|
\ge \frac{1}{4} \right\} \ge \frac{1}{2} \hspace{2cm} for \hspace{1cm}
k=1,...,n. \ter

\end{lemma} \hz

$Proof.$ First we observe that

\cen \inf \left\{ |<dM_k^{\infty } (\omega _1 ,...,\omega _k ),e_i >| :
i \in I(\omega _1 ,..., \omega _k ,- \omega _k ) \right\} = 2^{k-n-1}
\min (2^{n-k-1} +1, 2^{n-k} - 2^{n-k-1} ) \ge \frac{1}{4} \ter

for $1 \le k < n$. Then we use the fact that

\cen \# \left( \cap _1^n supp \hspace{0.2cm} dM_k^{\infty }(\omega ) \right)=1
\hspace{1cm} \mbox{for all} \hspace{1cm} \omega \in \{ -1,1 \} ^n \ter

to define $e(\omega )$ as the i-th unit vector, in the case if

\cen \{ i \} = \cap _1^n supp \hspace{0.2cm} dM_k^{\infty }(\omega ) \subseteq
\cap _2^n I(\omega _1 ,..., \omega _{k-1} ). \ter

For $1 \le k \le n-2$ we obtain \hz

$ \mu_n \left\{ \omega : |<dM_k^{\infty }(\omega ),e(\omega )>| \ge
\frac{1}{4} \right\} $ \hz

$ \pla \ge \mu_n \left\{ (\omega _1 ,..., \omega _n ) : |<dM_k^{\infty }
(\omega _1 ,..., \omega _k ),e(\omega _1 ,..., \omega _n )>| \ge \frac{1}{4} ,
\hspace{0.3cm} \omega _{k+1} = -\omega _k \right\} $ \hz

$ \pla \ge \mu_n \left\{ (\omega _1 ,..., \omega _n ) :
\inf \left\{ |<dM_k^{\infty }(\omega _1 ,..., \omega _k ),e_i>| : i \in
I(\omega _1 ,..., \omega _{k+1} ) \right\} \ge \frac{1}{4} , \hspace{0.3cm}
\omega _{k+1} = -\omega _k \right\} $ \hz

$ \pla = \mu_n \{ \omega _{k+1} = - \omega _k \} = \frac{1}{2} $. \hz

Since \pl $|<dM_k^{\infty }(\omega ),e(\omega )>| \ge \frac{1}{4} $ \pl for all
$\omega$ in the cases $k=n-1$ and $k=n$ the proof is complete. $\Box $ \hz

We deduce \hz

\begin{lemma} Let $n\ge 1$ be fixed. Then \\
(1) $ \| M_k^{\infty } (\omega )\| =1$ and $\| dM_l^{\infty } (\omega )\| =
    \frac{1}{2}$ \pll for $k=0,...,n$,\hspace{0.2cm} $l=1,...,n$, and all 
    $\omega \in \omet$,\\
(2) $\mu_n \left\{ \omega :\ew_{\vare }
    \| \sum_1^n \vare_k dM_k^{\infty } (\omega ) \|
    \ge \alpha n^{1/2} \right\} >\beta$ \pll for some $\alpha ,\beta >0$
    independent from $n$.
\end{lemma} \hz

$Remark.$ An inequality 
$\ew_{\vare }\| \sum_1^n \vare_k dM_k^{\infty } (\omega ) \| \ge \alpha n^{1/2}$
can not hold for all $\omega \in \omet$ since, for example,
\cen \| \sum_1^n \vare_k dM_k^{\infty} (1,1,...,1)\|  \le 
     \| \sum_1^n dM_k^{\infty} (1,1,...,1)\|          \le 
     \| \sigma_n \| \| \sum_1^n dM_k^1 (1,1,...,1)\| \le 2 \ter. 

$Proof$ of Lemma 4.2. Assertion (1) is trivial. We prove (2). For $t>0$ we
consider

\begin{eqnarray*}
\mu_n \left\{ \omega : \ew _{\vare } \| \sum_1^n \vare_k dM_k^{\infty }
(\omega ) \|_{l_{\infty }^{2^n }} > t n^{1/2} \right\} & \ge &
\mu_n \left\{ \omega : \| \ew_{\vare } | \sum_1^n \vare_k
dM_k^{\infty } (\omega ) |\|_{l_{\infty }^{2^n }} > t n^{1/2} \right\} \\
& \ge & \mu_n \left\{ \omega : \frac{1}{c_o} \| \left( \sum_1^n |
dM_k^{\infty } (\omega )|^2 \right) ^{1/2} \|_{l_{\infty }^{2^n }} > t n^{1/2}
\right\} \\
& = & \mu_n \left\{ \omega : \| \sum_1^n | dM_k^{\infty }
(\omega )|^2 \|_{l_{\infty }^{2^n }} > c_o^2 t^2 n \right\}.
\end{eqnarray*}

Denoting the last mentioned expression by $p_t$ the previous lemma yields

\begin{eqnarray*}
      p_t n + (1 - p_t) c_o^2 t^2 n 
&\ge& \ew_{\omega } \| \sum_1^n |dM_k^{\infty }(\omega ) |^2 \|          \\
&\ge& \ew_{\omega } <\sum_1^n |dM_k^{\infty }(\omega ) |^2 ,e(\omega )>  \\
& = & \ew_{\omega } \sum_1^n |<dM_k^{\infty }(\omega ),e(\omega )>|^2     \\
&\ge& \sum_1^n \frac{1}{16} \mu \left \{ \omega :
      |<dM_k^{\infty }(\omega ),e(\omega )>|^2 \ge \frac{1}{16} \right \} \\
&\ge& \frac{n}{32} 
\end{eqnarray*}

such that $p_t \ge \frac{1/32 - c_o^2 t^2}{1 - c_o^2 t^2}$ for
$c_o^2 t^2 < 1$. $\Box$ \hz

Lemmas 3.2 and 4.2 imply \hz

\begin{theorem} There are $\alpha ,\beta >0$ such that for all
non-superreflexive Banach spaces $X$, for all $\theta >0$, and for all
$n=1,2,...$ there exists a Walsh-Paley martingale $\{ M_k\}_0^n$ with values
in $B_X$,
\cen \inf_{1\le k\le n} \inf_{\omega } \| dM_k (\omega )\| \ge
\frac{1}{2(1+\theta )}, \pll \mbox{and} \pll
\mu_n \left\{ \omega :\ew_{\vare } \| \sum_1^n \vare_k dM_k(\omega ) \|
\ge \alpha n^{1/2} \right\} >\beta.\ter
\end{theorem} \hz

$Proof.$ We choose factorizations $\sigma_n = B_n A_n$ with
$A_n :l_1^{2^n} \rightarrow X$ , $B_n :X \rightarrow l_{\infty }^{2^n}$,
$\| A_n \| \le 1$ and $\| B_n \| \le 1+\min (1,\theta )$
(see \cite{J1}(Thm.4)). Defining
$M_k(\omega ):= A_n M_k^1(\omega )\in X$ we obtain
$\sup_{0\le k\le n,\omega \in \omet } \| M_k (\omega ) \| \le 1$ from Lemma 3.3
as well as
$\inf_{1\le k\le n,\omega \in \omet } \| dM_k (\omega ) \| \ge
\frac{1}{2(1+\theta )}$ and

\begin{eqnarray*}
\mu_n \left\{ \omega : \ew _{\vare } \| \sum_1^n \vare_k dM_k(\omega )
\|_X > \frac{\alpha }{2} n^{1/2} \right\}
& \ge & \mu_n \left\{ \omega : \ew _{\vare } \| \sum_1^n \vare_k
dM_k(\omega )\|_X > \frac{\alpha }{\| B_n \| } n^{1/2} \right\} \\
& \ge & \mu_n \left\{ \omega : \ew _{\vare } \| \sum_1^n \vare_k dM_k^{\infty }
(\omega )\| _{l_{\infty }^{2^n}} > \alpha n^{1/2} \right\} \ge \beta
\end{eqnarray*}

according to Lemma 4.2.$\Box$ \hz

For Banach spaces of type 2 we get \hz

\begin{theorem} For any Banach space $X$ of type 2 the following assertions 
are equivalent.\\
(1) $X$ is not superreflexive.\\
(2) $\frac{1}{c} n^{1/2} \le \rue X) \le c n^{1/2}$ \pll for $n=1,2,...$
    and some $c>0$.\\
(3) $\frac{1}{c'} n^{1/2} \le \rue X)$ \pll for $n=1,2,...$ and some
    $c'>0$.\\
\end{theorem} \hz
 
$Proof.(1) \Rightarrow (2)$ follows from Theorem 4.3 and Lemma 3.4. 
$(3) \Rightarrow (1)$. We assume $X$ to be superreflexive and find
(\cite{J2}, cf.\cite{P1}(Thm1.2,Prop.1.2)) $\gamma >0$ and $2\le s<\infty$
such that

\cen \left( \sum_{k\ge 0} \| dM_k\|_{L_2^X }^s \right)^{1/s} \le \gamma
\sup_k \| M_k\|_{L_2^X } \ter

for all martingales in $X$. This martingale cotype implies

\begin{eqnarray*}
\left( \ew_{\vare ,\omega } \| \sum_1^n \vare_k dM_k (\omega )\|^2
\right)^{1/2}
& \le & T_2(X) \left( \ew_{\omega } \sum_1^n \| dM_k (\omega )\|^2
\right)^{1/2} \\
& \le & T_2(X) n^{1/2 - 1/s} \gamma \| M_n - M_0\|_{L_2^X}
\end{eqnarray*}

which contradicts $\rue X) \ge \frac{1}{c'} n^{1/2}$.$\Box$ \hz

$Remark.$ Corollary 5.4 will demonstrate that the asymptotic
behaviour of $\rue X)$ can not characterize the superreflexivity of $X$
in the case that $X$ is of type $p$ with $p<2$. Namely, according to
Theorem 5.4 for all $1<p<2<q<\infty$ there is a superreflexive Banach
space $X$ of type $p$ and of cotype $q$ with 
$\rue X) \asymp n^{\frac{1}{p} -\frac{1}{q}}$. On the
other hand, if $\frac{1}{p} -\frac{1}{q} \ge \frac{1}{2}$ then we can find
a non-superreflexive Banach space $Y$ such that 
$\rue Y) \asymp n^{\frac{1}{p} -\frac{1}{q}}$ (add a non-superreflexive
Banach space of type 2 to X). \hz

Finally, we deduce the random unconditional constants of the summation
operators $\sigma_n$, $\sigma$, and $\Phi$ defined in the beginning of this
section. To this end we need the type 2 property of these operators. From
\cite{J1} and \cite{J3} or \cite{P-X} as well as \cite{X} we know the much 
stronger results, that $\sigma$ and the usual summation operator from 
$L_1 [0,1]$ into $L_{\infty } [0,1]$ can be factorized through a type 2 space.
We want to present a very simple argument for the type 2 property
of the operator $\Phi$ which can be extended to some other "integral 
operators" from $C[0,1]'$ into $l_{\infty } ([0,1])$.\hz

\begin{lemma} The operator $\Phi :C[0,1]' \rightarrow l_{\infty } ([0,1])$ is 
of type 2 with $T_2(\Phi ) \le 2$.
\end{lemma}

$Proof.$ First we deduce the type 2 inequality for Dirac-measures. Let
$\lambda_1 ,...,\lambda_n \in \kz$, $t_1 ,...,t_n \in [0,1]$, whereas we
assume $0\le t_{k_1} = ... = t_{l_1} < t_{k_2} = ... = t_{l_2} < ... <
t_{k_M} = ... = t_{l_M} \le 1$, and let $\delta_{t_1} ,...,\delta_{t_n}$
the corresponding Dirac-measures. Then, using Doob's inequality, we obtain

\begin{eqnarray*}
\kla \ew_{\vare } \| \sum_1^n \Phi \vare_j \lambda_j \delta_{t_j} \|^2 
\mer^{1/2}
& = & \kla \ew_{\vare } \sup_t | \sum_{i=1}^M \kla \sum_{k_i}^{l_i} \vare_j
\lambda_j \mer \delta_{t_{k_j}} ([0,t]) |^2 \mer^{1/2} \\ 
& = & \kla \ew_{\vare } \sup_{1\le m\le M} | \sum_{i=1}^m \kla 
\sum_{k_i}^{l_i} \vare_j \lambda_j \mer |^2 \mer^{1/2} \\
& \le & 2 \kla \ew_{\vare } | \sum_{i=1}^M \kla \sum_{k_i}^{l_i} \vare_j 
\lambda_j \mer |^2 \mer^{1/2} \\
& = & 2 \kla \sum_1^n |\lambda_j |^2 \mer^{1/2}.
\end{eqnarray*}

Hence

\cen \kla \ew_{\vare } \| \sum_1^n \Phi \vare_j \lambda_j \delta_{t_j} \|^2 
\mer^{1/2} \le 2 \kla \sum_1^n \| \lambda_j \delta_{t_j} \|^2 \mer^{1/2}. \ter

In the next step for any $\mu \in C[0,1]'$ we find a sequence of point measures
(finite sums of Dirac-measures) $\{ \mu^m \}_{m=1}^{\infty} \subset C[0,1]'$ 
such that
$\sup_m \| \mu^m \| \le \| \mu \|$ and $\lim_m \mu^m ([0,t]) = \mu ([0,t])$
for all $t\in [0,1]$ (take, for example,
$\mu^m := \sum_{i=1}^{2^m} \delta_{\frac{i-1}{2^m}} \mu (I_i^m )$
with $I_1^m := [0,\frac{1}{2^m} ]$ and 
$I_i^m := (\frac{i-1}{2^m} ,\frac{i}{2^m} ]$ for $i>1$). Now, assuming
$\mu_1 ,...,\mu_n \in C[0,1]'$ we choose for each $\mu_j$ a sequence
$\{ \mu_j^m \}_{m=1}^{\infty}$ of point measures in the above way and obtain
 
\begin{eqnarray*}
\kla \ew_{\vare } \| \sum_1^n \Phi \vare_j \mu_j \|^2 \mer^{1/2}
& = & \kla \ew_{\vare } \sup_t |\sum_1^n \vare_j \mu_j ([0,t]) |^2 \mer^{1/2}\\
& = & \kla \ew_{\vare } \sup_t \lim_m |\sum_1^n \vare_j \mu_j^m ([0,t]) |^2 
\mer^{1/2}\\
& \le & \limsup_m \kla \ew_{\vare } \sup_t |\sum_1^n \vare_j \mu_j^m ([0,t]) 
|^2 \mer^{1/2}.
\end{eqnarray*}

Using the type 2 inequality for Dirac measures and an extreme point argument
we may continue to

\cen \kla \ew_{\vare } \| \sum_1^n \Phi \vare_j \mu_j \|^2 \mer^{1/2} \le
2 \limsup_m \kla \sum_1^n \| \mu_j^m \|^2 \mer^{1/2} \le
2 \kla \sum_1^n \| \mu_j \|^2 \mer^{1/2}. \Box \ter \hz

As a consequence we obtain \hz

\begin{theorem} There is an absolute constant $c>0$ such that for all
$n=1,2,...$
\cen \frac{1}{c} n^{1/2} \le \rue \sigma_n ) \le \rue \sigma ) \le
\rue \Phi ) \le c n^{1/2}. \ter
\end{theorem} \hz

$Proof. \frac{1}{c} n^{1/2} \le \rue \sigma_n )$ is a consequence of
Lemma 3.3 and 4.2. $\rue \sigma_n ) \le \rue \sigma ) \le \rue \Phi )$
is trivial. Finally, Lemma 4.5 and Lemma 3.4 imply $\rue \Phi ) \le 
4 n^{1/2} . \Box$ \hz

\begin{cor} There is an absolute constant $c>0$ such that for all
$n=1,2,...$
\cen \frac{1}{c} n^{1/2} \le 
\kla \ew_{\vare, \omega} \| \sum_1^n \vare_k dM_k^{\infty } (\omega ) \|^2
\mer^{1/2} = 
\kla \ew_{\vare, \omega} \| \sum_1^n \vare_k |dM_k^{\infty } (\omega )| 
\hspace{.2cm} \|^2 \mer^{1/2} \le c n^{1/2}. \ter
\end{cor} \hz

$Proof.$ This immediately follows from Lemma 4.2, Theorem 4.6, and
$dM_k^{\infty } (\omega )= \omega_k |dM_k^{\infty } (\omega )|. \Box$ \hz

\section{An example}
\setcounter{lemma}{0}

We consider an example of Bourgain to demonstrate that for all
$0\le \alpha <1$ there is a superreflexive Banach space $X$ with 
$\rue X) \asymp n^\alpha$. Moreover, the general principle of this
construction allows us to show that $RU\!M\!D_n^1 (\kz ) \asymp n$ mentioned
in section 2 of this paper. \\ 
The definitions concerning upper p- and lower-q estimates of a
Banach space as well as the modulus of convexity and smoothness,
which we will use here, can be found in \cite{L-T}.\hz

Let us start with a Banach space $X$ and let us consider the function space
$X_{\omet } := \{ f:\omet \rightarrow X \}$ equipped with some norm
$\| \pl \| = \| \pl \|_{X_{\omet }}$. For a fixed $f\in X_{\omet }$ we define

\cen M^f : \omet \rightarrow X_{\omet } \hspace{1cm} \mbox{by} \hspace{1cm}
M^f (\omega ) := f_{\omega } \ter

where $f_{\omega } (\omega '):=f(\omega \omega ')$ \hspace{.3cm}
($\omega \omega ' := (\omega_1 \omega_1 ' ,...,\omega_n \omega_n ' )$ for
$\omega =(\omega_1 ,...,\omega_n )$ and $\omega ' =(\omega_1 ' ,...,
\omega_n ' )$). Setting $M_k^f := \ew (M_n^f |{\cal F} _k )$ we obtain a 
martingale $\{ M_k^f \}_{k=o}^n$ with values in $X_{\omet }$ generated by the 
function $f \in X_{\omet }$. Furthermore, putting $f_n := f$,

\cen f_k (\omega ) := \ew (f | {\cal F} _k )(\omega ) = \frac{1}{2^{n-k}}
\sum_{\omega_{k+1} ' = \pm 1} ... \sum_{\omega_n ' = \pm 1}
f_n (\omega_1 ,...,\omega_k ,\omega_{k+1} ',...,\omega_n '), \ter

$df_k := f_k - f_{k-1}$ for $k\ge 1$, and $df_0 = f_0$, it yields

\cen \kla \sum_0^n \alpha_k df_k \mer _{\omega } =
\sum_0^n \alpha_k dM_k^f (\omega ) \hspace{0.5cm} \mbox{for all} \hspace{0.5cm}
\omega \in \omet \hspace{0.5cm} \mbox{and all} \hspace{0.5cm}
\alpha_0 ,..., \alpha_n \in \kz .\ter

The following lemma is now evident.\hz

\begin{lemma} Let $f\in X_{\omet }$ and let $\{ M_k^f \}_0^n$ be the
corresponding martingale. If $\| \pl \| = \| \pl \|_{X_{\omet }}$ is 
translation invariant then
$\| \sum_0^n \alpha_k dM_k^f (\omega ) \| = \| \sum_0^n \alpha_k df_k \|$
for all $\omega \in \omet$ and all $\alpha_0 ,....,\alpha_n \in \kz$.
\end{lemma} \hz

First we deduce \hz

\begin{cor} There exists $c>0$ such that $\frac{n}{c} \le \ruei \kz ) \le 
\ruei X) \le cn$ for all $n=1,2,...$ and all Banach spaces $X$.
\end{cor} \hz

$Proof.$ We consider $\kz_{\omet }$ with $\| f\| := \sum_{\omega } |f(\omega )|$
such that $\kz_{\omet } = l_1 (\omet )$. Defining $f \in l_1 (\omet )$ as
$f := \chi _{\{ (1,...,1)\} }$ it follows that $f_{\omega } =
\chi_{\{ \omega \} }$. It is clear that the isometry $I:l_1 (\omet ) \rightarrow
l_1^{2^n}$ with $If_{\omega } := e_{i(\omega )}$ ($e_{i(\omega )}$ is the
$i(\omega )$-th unit vector where $i(\omega )$ is defined as in section 3 of
this paper) transforms the martingale $\{ M_k^f\}$ into the martingale
$\{ M_k^1\}$ from section 3 by $IM_k^f(\omega ) = M_k^1(\omega )$ for all
$\omega \in \omet$. Combining Lemma 5.1 and Lemma 3.3 yields

\cen \inf_{\omega } \ew_{\vare } \| \sum_1^n \vare_k df_k \|_{l_1 (\omet )}
\ge {\alpha}n \hspace{1cm} \mbox{and} \hspace{1cm} \|
f-f_0\|_{l_1 (\omet )} \le 2.
\ter

Consequently, $\ruei X) \ge \ruei \kz ) \ge \frac{\alpha}{2} n$. On the other
hand we have $\ruei X) \le 2n$ in general.$\Box$ \hz

Now, we treat Bourgain's example \cite{Bou} .\hz

\begin{theorem} For all $1<p<q<\infty$ and $n \in \nz$ there exists a
function lattice $X_{pq}^{2n} = \kz_{\omett }$ such that\\
(1) $X_{pq}^{2n}$ has an upper $p$- and a lower $q$-estimate with the
    constant 1,\\
(2) there exists a Walsh-Paley martingale $\{ M_k\}_0^{2n}$ with values
    in $B_{X_{pq}^{2n}}$ and
    \cen \inf_{\omega } \ew_{\vare } \| \sum_1^{2n} \vare_k dM_k (\omega )\|
    \ge c(2n)^{\frac{1}{p} -\frac{1}{q} } \ter
    where $c>0$ is an absolute constant independent from $p$,$q$, and $n$.
\end{theorem} \hz

$Proof.$ In \cite{Bou} (Lemma 3) it is shown that there is a lattice norm
$\| \pl \| = \| \pl \|_{\kz_{\omett }}$ on $\kz_{\omett }$ which satisfies (1),
such that there exists a function $\phi \in \kz_{\omett }$ with

\cen \| \phi \| \le \vare^{1-\frac{p}{q} } \hspace{0.5cm}
(\vare = 2n^{-\frac{1}{p} }) \hspace{1cm} \mbox{and} \hspace{1cm}
\| \kla \sum_0^{2n} |d\phi_k |^2 \mer^{1/2} \| \ge \frac{1}{2} \ter

\cite{Bou} (Lemma 4 and remarks below, $\vare = 2n^{-1/p}$ is taken from
the proof of Lemma 4). Since $\| \pl \|$ is translation invariant Lemma 5.1
implies

\cen \| M_k^{\phi } (\omega ) \| = \| M_k^{\phi }\|_{L_2^{X_{pq}^{2n}}} \le 
\| M_n^{\phi }\|_{L_2^{X_{pq}^{2n}}} = \| \phi \|
\le 4 (2n)^{\frac{1}{q} -\frac{1}{p}} \ter

and

\cen \ew_{\vare } \| \sum_0^{2n} \vare_k dM_k^{\phi } (\omega ) \| =
\ew_{\vare } \| \sum_0^{2n} \vare_k d{\phi }_k \| \ge
\| \ew_{\vare } | \sum_0^{2n} \vare_k d{\phi }_k | \hspace{0.2cm} \| \ge
\frac{1}{A} \| \kla \sum_0^{2n} |d{\phi }_k|^2 \mer^{1/2} \| \ge \frac{1}{2A}.
\Box \ter \hz

As usual, in the following the phrase "the modulus of convexity (smoothness) of
$X$ is of power type $r$" stands for "there is some equivalent norm on $X$ 
with the modulus of convexity (smoothness) of power type $r$". Now, similarly 
to \cite{Bou} we apply a standard procedure to the above finite-dimensional 
result. 

\begin{cor} (1) For all $1<p<2<q<\infty$ there is a Banach space $X$  
with the modulus of convexity of power type $q$ and the modulus of smoothness 
of power type $p$, and a constant $c>0$ such that
\cen \frac{1}{c} n^{\frac{1}{p} - \frac{1}{q}} \le \rue X) \le c 
n^{\frac{1}{p} - \frac{1}{q}} \hspace{.5cm} \mbox{for} \hspace{0.5cm}
n=1,2,...\ter
(2) There is a Banach space $X$ with the modulus of convexity of power type 
$q$ and the modulus of smoothness of power type $p$ for all $1<p<2<q<\infty$, 
and $\rue X) \rightarrow_{n \rightarrow \infty} \infty$.
\end{cor}\hz

$Proof.$ For sequences $P=\{ p_n \}$ and $Q=\{ q_n \}$ with

\cen 1<p_1 \le p_2 \le ... \le p_n \le ... <2<... \le q_n \le ... \le q_2
\le q_1 < \infty \ter

we set $X_{PQ} := \bigoplus_2 X_{p_n q_n}^{2n}$ and obtain that
$X_{PQ}$ satisfies an upper $p_k$- and a lower $q_k$-estimate for all
$k$. According to a result of Figiel and Johnson (cf. \cite{L-T} (II.1.f.10))
$X_{PQ}$ has the modulus of convexity of power type $q_k$ and the modulus
of smoothness of power type $p_k$ for all $k=1,2,...$. Furthermore,
\cite{P1}(Theorem 2.2) implies

\begin{eqnarray*}
\sup_{\vare_1 \pm 1,...,\vare_n \pm 1} \| \sum_1^n \vare_l dM_l \|_{L_2^X} 
& \le & c_{p_k} \kla \sum_1^n \| dM_l \|_{L_2^X}^{p_k} \mer^{1/{p_k}} \\
& \le & c_{p_k} n^{\frac{1}{p_k} - \frac{1}{q_k}}
\kla \sum_1^n \| dM_l \|_{L_2^X}^{q_k} \mer^{1/{q_k}} \\
& \le & c_{p_k} d_{q_k} n^{\frac{1}{p_k} - \frac{1}{q_k}}
\| \sum_1^n dM_l \|_{L_2^X}
\end{eqnarray*}

for all martingales $\{ M_l\}$ with values in $X_{PQ}$ such that
$\rue X_{PQ}) \le c_{p_k} d_{q_k} n^{\frac{1}{p_k} - \frac{1}{q_k}}$.
On the other hand, from Theorem 5.3 we obtain

\cen c(2n)^{\frac{1}{p_n} -\frac{1}{q_n}} \le \ruee X_{p_n q_n}^{2n} )\le
\ruee X_{PQ} ). \ter

Now, setting $p_k \equiv p$ and $q_k \equiv q$ we obtain (1). Choosing the
sequences in the way that $p_k \rightarrow_{k \rightarrow \infty} 2$, 
$q_k \rightarrow_{k \rightarrow \infty} 2$, and
$n^{\frac{1}{p_n} - \frac{1}{q_n}} \rightarrow_{n \rightarrow \infty} \infty$ 
assertion (2) follows.$\Box$


\begin{thebibliography}{99}

 \bibitem{Al} D.J.Aldous , Unconditional bases and martingales in $L_p(F)$,
 Math.Proc.Camb.Phil.Soc.85(1979), 117-123.

 \bibitem{Bou} J.Bourgain , Some remarks on Banach spaces in which
 martingale difference sequences are unconditional, Ark.Mat. 21(1983), 163-168.

 \bibitem{Bu} D.L.Burkholder , Martingales and Fourier analysis in Banach
 spaces, Probability and analysis (Varena, 1985),LNM 1206(1986), 61-108.

 \bibitem{El} J.Elton , Sign-embeddings on $l_1^n$,
 Trans.Amer.Math.Soc.279(1983), 113-124.

 \bibitem{Ga} D.J.H.Garling , Random martingale transform inequalities,
 Probability in Banach spaces 6. Proceedings of the Sixth international
 conference, Sandbjerg, Denmark 1986, Birkhaeuser 1990.

 \bibitem{G} S.Geiss , Antisymmetric tensor products of absolutely p-summing
 operators, to appear in J.Appr.Th.

 \bibitem{J1} R.C.James , Super-reflexive Banach spaces, Can.J.Math.
 5(1972), 896-904.

 \bibitem{J2} R.C.James , Super-reflexive spaces with bases,
 Pacific J.Math.41(1972), 409-419. 
 
 \bibitem{J3} R.C.James , Non-reflexive spaces of type 2, Isr.J.Math.
 30(1978), 1-13.
 
 \bibitem{L-T} J.Lindenstrauss and L.Tzafriri , Classical Banach spaces I/II,
 Springer,New York-Berlin-Heidelberg 1977/79.

 \bibitem{M} B.Maurey , Systeme de Haar. Sem. Maurey-Schwartz 1974-1975,
 Ecole Polytechnique, Paris.

 \bibitem{M-P} B.Maurey and G.Pisier , Series de variables aleatoires 
 vectorielles independantes et proprietes geometriques des espaces de Banach, 
 Stud. Math. 58(1976), 45-90.

 \bibitem{Pa} A.Pajor , Prolongement de $l_1^n$ dans les espaces de Banach
 complexes, C.R.Acad.Sci.Paris Ser. I Math. 296(1983), 741-743.

 \bibitem{P1} G.Pisier , Martingales with values in uniformly convex spaces,
 Isr.J.Math. 20(1975), 326-350.

 \bibitem{P2} G.Pisier , Un example concernant la super-reflexivite,
 Sem. Maurey-Schwartz 1974-1975, Ecole Polytechnique, Paris.

 \bibitem{P-X} G.Pisier,Q.Xu , Random series in the real interpolation
 spaces between the spaces $v_p$, Geometrical aspects of functional
 analysis (1985/86), LNM 1267(1987), 185-209.

 \bibitem{T-J} N.Tomczak-Jaegermann , Banach-Mazur distances and finite
 dimensional operator ideals, Pitman 1988.

 \bibitem{X} Q.Xu , Espaces d'interpolation reels entre les espaces $V_p$:
 Proprietes geometriques et applications probabilistes, Theses, Universite
 Paris 6, 1988.

\end{thebibliography}
\end{document}